\documentclass[a4paper, 12pt]{article}
\usepackage{amssymb}

\usepackage{amsmath,amsthm}
\date{}

\title{\LARGE \bf
Algebraic Characterizations of Consensus Problems for Networked
Dynamic Systems
\thanks{This work was supported by
the National Natural Science Foundation of China (No.  10372002
and No.  60274001) and the National Key Basic Research and
Development Program (No.  2002CB312200).}}

\author{Long Wang,  Hong Shi, Feng Xiao, and Aiping Wang\\
Intelligent Control Laboratory, Center for Systems and Control,\\
  Department of Mechanics and
Engineering Science,\\        Peking University,  Beijing 100871,
China.\\     E-mail address: {\tt\small longwang@pku.edu.cn} }

\newtheorem{dl}{Theorem}[section]
\newtheorem{tl}{Corollary}[section]
\newtheorem{yl}{Lemma}[section]

\newtheorem{rmrk}{Remark}[section]
\newtheorem{ff}{Method}
\newtheorem{lz}{Example}[section]

\newcommand{\R}{\mathbf{R}}
\newcommand{\n}{\bar{n}}

\begin{document}

\maketitle \thispagestyle{empty} \pagestyle{empty}

%%%%%%%%%%%%%%%%%%%%%%%%%%%%%%%%%%%%%%%%%%%%%%%%%%%%%%%%%%%%%%%%%%%%%%%%%%%%%%%%
\begin{abstract}

In this paper,  we study the consensus problem for networked
dynamic systems with arbitrary initial states,  and present some
structural characterization and direct construction of consensus
functions.  For the consensus problem under similar
transformation, we establish some necessary and sufficient
conditions by exploiting the structure of consensus functions.
Finally,  we discuss the consensus problem for dynamic systems
under switching by using the common Lyapunov function method.

\end{abstract}


%%%%%%%%%%%%%%%%%%%%%%%%%%%%%%%%%%%%%%%%%%%%%%%%%%%%%%%%%%%%%%%%%%%%%%%%%%%%%%%%

\section{INTRODUCTION}

In recent years,  the study of synchronization and coordination of
multi-agent systems has attracted many researchers.  It has broad
applications in  cooperative control of unmanned aerial vehicles,
scheduling of automated highway systems,  formation control of
satellite clusters,  and distributed optimization of multiple
mobile robotic systems(\cite{a1}--\cite{R.  Bachmayer and N.  E.
Leonard}).

Consensus problem has a long history. On many occasions,  a group
of dynamic agents in multi-agent/multi-robot systems need to reach
an agreement on certain quantities of interest.  For example,
flock of birds tends to synchronize in migration in order to
resist external aggression and reach their destination.  Robots
need to arrive at agreement so as to accomplish some complicated
tasks. Investigation of such problems is of significance  in
theory and in practice.

Consensus problem was introduced and formally stated by
\cite{c1}-\cite{c2}.  In \cite{c1}, the basic definitions were
given and average consensus problem was studied for networks with
both switching topology and time delays.  In this paper,  we
generalize the consensus problem and formulate it in a more
general form. \cite{c1} tackled this problem mostly by graph
theory and assumed that the state of each agent is a real scalar.
However, in most cases, the quantities of each agent are very
complex and many aspects should be considered.  For example,  the
quantities might be position,  velocity, temperature,  momentum,
voltage, mass, energy and so on. Furthermore,  these quantities
might not be independent. Hence, it is natural to extend the
domain $\mathbf{R}$ of the state of each agent to $\mathbf{R}^m$.
Therefore, all the original definitions for consensus problem
should be modified correspondingly. In this paper, this kind of
consensus problem is studied by using linear algebra theory
(\cite{b1, R. Horn and C. R. Johnson}) as basic tool, and some
interesting structural characterizations are established.

This paper is organized as follows: In section II, we generalize
consensus problem  and establish some necessary and sufficient
conditions for a linear dynamic system that solves a consensus
problem with arbitrary initial state.  In section III, we focus on
the structural characterizations of consensus functions and
present a simple and constructive method to obtain  consensus
functions.  Furthermore,  a necessary and sufficient condition for
a dynamic system that solves the average consensus problem with
arbitrary initial state is given.  In Section IV, the consensus
problem under similar transformation is discussed.  In Section V,
the systems that solve a consensus problem under arbitrary
switching(\cite{a3, b2, b3}) are characterized. Finally, we
summarize our main contribution in Section VI. For convenience,
some concepts and results in graph theory are given in the
Appendix.

%%%%%%%%%%%%%%%%%%%%%%%%%%%%%%%%%%%%%%%%%%%%%%%%%%%%%%%%%%%%%%%%%%%%%%%%%%%%%%%%

\section{PRELIMINARIES}
In order to introduce the generalized consensus concept,  we
consider the following linear dynamic system:
\begin{equation}\label{xt1}
    \left[%
\begin{array}{c}
  \dot{x}_1 \\
  \dot{x}_2 \\
  \vdots \\
  \dot{x}_n \\
\end{array}%
\right]
=\left[%
\begin{array}{cccc}
  A_{11} & A_{12} & \cdots & A_{1n} \\
  A_{21} & A_{22} & \cdots & A_{2n} \\
  \vdots & \vdots & \ddots & \vdots \\
  A_{n1} & A_{n2} & \cdots & A_{nn} \\
\end{array}%
\right]
\left[%
\begin{array}{c}
  {x_1} \\
  {x_2} \\
  \vdots \\
  {x_n} \\
\end{array}%
\right],
\end{equation}
where $x_i\in \mathbf{R}^{m}$,  $A_{ij} \in \mathbf{R}^{m\times
m}$, $i, j=1, 2, \cdots, n$.  System (\ref{xt1}) can be written in
the following form
\begin{equation}\label{xt11}
    \dot{x}= Ax,
\end{equation}
where $A=[A_{ij}]$ and $ x=[x_1^T, x_2^T, \cdots, x_n^T]^T$.

We say $x_i$ and $x_j$ {\it agree} if and only if
$x_i=x_j$(component-wise).  Let $\n=\{1,  2,  \cdots,  n\}$.  We
say system (\ref{xt1}) has reached a {\it consensus} if and only
if $x_i=x_j $ for all $i \neq j, i, j\in \n$. The common value of
$x_i$ ($i=1, 2, \cdots, n$) is called the {\it group decision
value}. Let $\chi: \R^{mn} \to \R^m$ be a function of $n$ vectors
$x_1, x_2, \cdots, x_n (x_i \in \R^m)$ and $x(0)$ denote the
initial state of the system. We say dynamic system {\it solves the
$\chi$-consensus problem} if and only if there exists an
asymptotically stable equilibrium $x^*=[x^{*T}_1,  \cdots,
x^{*T}_n]^T$ of system (\ref{xt1}) satisfying $x_i^* = \chi(x(0))
\in\R^m$ for all $i \in \n$. The function $\chi$ is called {\it
consensus function}. The special cases of $\chi(x)={\mathrm
{Ave}}(x)= (\sum\nolimits_{i=1}^n x_i)/n$,
$\chi(x)=\mathop{\max}\limits_{\|x_i\|}(x_i)$,
$\chi(x)=\mathop{\min}\limits_{\|x_i\|}(x_i)$ are called {\it
average-consensus},  {\it max-consensus} and {\it min-consensus},
respectively, due to their broad applications in distributed
decision making for multi-agent systems.  If we have $x_i^*=x^*_j$
for all $i\neq j, i,  j \in \n$,  and $x_i^*$ only relies on
initial state $x(0)$, we say that the system {\it solves a
consensus problem}.

Here, we are interested in the system $\dot{x}=Ax$ which solves
the $\chi$-consensus problem for some consensus function $\chi$
and for any $x(0) \in \R^{mn}$. For such systems, there are some
necessary properties. For example,  if $x \in N(A)$,  then
$x=\mathbf{1}\otimes b$, where $N(A)$ is the null space of $A$,
$\mathbf{1}=[1,  1, \cdots, 1]^T\in{R^n}$,  $\otimes$ is the
Kronecher product,  and $b \in \R^m$ is a constant vector.
Furthermore,  for any initial state $x(0)$, the solution of the
system converges asymptotically to some equilibrium. Denote the
range (column space) of $A$ by $R(A)$. We have the following
theorem.

\begin{dl}\label{dl21}
System (\ref{xt1}) solves a consensus problem for any initial
state $x(0)$ if and only if $R(A)=R(A^2)$  and each eigenvalue of
$A$ is $0$ or has negative real part.  Moreover, if $0$ is  an
eigenvalue of $A$, then for any $x\in N(A)$, there exists a vector
$b\in \R^m$ such that $x=\mathbf{1}\otimes b$.
\end{dl}

\begin{proof}{\it Necessity}.  If $R(A)\neq R(A^2)$,  then we have
rank$(A)> $ rank$(A^2)$.  Hence, there exists a vector
$y\in{\R^{mn}}$, $y\in N(A^2)$ but $y \notin N(A)$. The solution
of system (\ref{xt1}) with the initial value $y$ is $x=e^{At}y
=(I_{mn} + At + \frac{A^2t^2}{2} + \cdots )y=y +tAy$, where
$I_{mn}$ is identical matrix of order $mn$. Obviously, $x$ does
not converge to any equilibrium when $t\to\infty$,  which is a
contradiction. Therefore $R(A)=R(A^2)$.

{\it Sufficiency}.  If all eigenvalues of $A$ have negative real
parts, then the system is asymptotically stable and all solutions
converge to 0,  i.e.,  it solves a consensus problem. If $0$ is an
eigenvalue of $A$,  then there exists an invertible matrix $T$ by
$R(A)=R(A^2)$, such that
\begin{equation}\label{jor1}
    A=T^{-1}\left[%
\begin{array}{cccccc}
  0 &  &  &  &  &  \\
   & \ddots &  &  &  &  \\
   &  & 0 &  &  &  \\
   &  &  & J_2 &  &  \\
   &  & &  & \ddots &  \\
   &  &  &  &  & J_s \\
\end{array}%
\right]T,
\end{equation}
where $J_2,  \cdots,  J_s $ are Jordan blocks,  and the eigenvalue
of $J_i$ has negative real part.  So the system converges
asymptotically to some equilibrium,  and since for any $x\in
N(A)$, there exists a vector $b\in \R^m$ such that
$x=\mathbf{1}\otimes b$,  the system solves a consensus problem
for any initial state.
\end{proof}

\begin{tl}\label{tl22}
 If system (\ref{xt1}) solves a consensus problem for any initial state,  then dim$(N(A))\leq m. $
\end{tl}

\begin{tl}\label{tl23}
 If  system (\ref{xt1}) solves a consensus problem for any initial state,  then
 \[\R^{mn}=N(A)\oplus R(A),\]
 where $\oplus$ is the operator of direct sum,  and
 \[A(R(A))=R(A), \; \; A(N(A))=\{0\}.\]
\end{tl}
\begin{proof}
From Theorem \ref{dl21}, we obtain rank$(A)=$ rank$(A^2)$.  The
remained proof is trivial.
\end{proof}

In order to investigate system (\ref{xt1}) more insightfully,
(\ref{xt1}) can be formulated in the following form
\begin{equation}\label{xt12}
\begin{array}{rcl}
    \left[%
\begin{array}{c}
  \dot{x}_1 \\
  \dot{x}_2 \\
  \vdots \\
  \dot{x}_n \\
\end{array}%
\right]
&=&\left[%
\begin{array}{cccc}
  C_{11} &  & &  \\
   & C_{22} & &  \\
   &  & \ddots &  \\
  &  &  & C_{nn} \\
\end{array}%
\right]
\left[%
\begin{array}{c}
  {x_1} \\
  {x_2} \\
  \vdots \\
  {x_n} \\
\end{array}%
\right]\\
& +&
\left[%
\begin{array}{cccc}
  D_{11} & D_{12} & \cdots & D_{1n} \\
  D_{21} & D_{22} & \cdots & D_{2n} \\
  \vdots & \vdots & \ddots & \vdots \\
  D_{n1} & D_{n2} & \cdots & D_{nn} \\
\end{array}%
\right]
\left[%
\begin{array}{c}
  {x_1} \\
  {x_2} \\
  \vdots \\
  {x_n} \\
\end{array}%
\right],
\end{array}
\end{equation}
where $C_{ii},  D_{ij} \in \mathbf{R}^{m\times m}$  such that
$\sum\nolimits_{j=1}^nD_{ij}=0$ for all $i, j \in \n$.

\begin{dl}\label{dl24}
System (\ref{xt1}) solves a consensus problem for any initial
state $x(0)$ if and only if
\[{\mathrm{dim}}N(A)={\mathrm{dim}}N(A^2)={\mathrm{dim}}N([C_{11}^T, C_{22}^T, \cdots,
C_{nn}^T]^T)\] and each eigenvalue of $A$ is $0$ or has negative
real part.
\end{dl}

\begin{proof}
We only need to prove that the condition
\[{\mathrm {dim}}N(A)={\mathrm {dim}}N(A^2)={\mathrm {dim}}N([C_{11}^T, C_{22}^T, \cdots,
C_{nn}^T]^T)\] is equivalent to the conditions that $R(A)=R(A^2)$
and for any $x\in N(A)$, there exists a vector $b\in \R^m$ such
that $x=\mathbf{1}\otimes b$.

Obviously, dim$N(A)$=dim$N(A^2)$ is equivalent to $R(A)=R(A^2)$.

(a) Suppose that dim$N(A)=r$, and for any $x\in N(A)$, there
exists a vector $b\in \R^m$ such that $x=\mathbf{1}\otimes b$. The
equation $Ax=0$ must have $r$ linearly independent solutions,
which implies that there are $r$ linearly independent vectors
$b_1, b_2, \cdots,  b_r \in \R^m$ such that $A(\mathbf{1}\otimes
b_i)=0$ for any $i \in \{1, 2, \cdots, r\}$. Substituting
$\mathbf{1}\otimes b_i$ into (\ref{xt12}), we obtain $C_{ii}b_j=0$
for any $i \in \n, j\in \{1,  2,  \cdots, r\}$. Therefore
dim$N([C_{11}^T, C_{22}^T, \cdots, C_{nn}^T]^T)\geq r$. But if
dim$N([C_{11}^T, C_{22}^T, \cdots, C_{nn}^T]^T)>r$,  then the
number of linearly independent solutions of  the equation
$[C_{11}^T, C_{22}^T, \cdots, C_{nn}^T]^Tx=0$ is more than $r$,
which implies that dim$(N(A))>r$, which contradicts our
assumption. Therefore dim$N([C_{11}^T, C_{22}^T, \cdots,
C_{nn}^T]^T)=r$.

(b) If dim$N([C_{11}^T, C_{22}^T, \cdots, C_{nn}^T]^T)=r$,  then
there are $r$ linearly independent solutions $b_1,  b_2,  \cdots,
b_r \in \R^m$ of equation $[C_{11}^T, C_{22}^T, \cdots,
C_{nn}^T]^Tx=0$. Thus the equation $Ax=0$   has $r$ independent
solutions $\mathbf{1}\otimes b_1$,  $\mathbf{1}\otimes b_2$,
$\cdots$, $\mathbf{1}\otimes b_r$.  Since dim$N(A)=r$,  we obtain
that for any $x\in N(A)$, there exists a vector $b\in \R^m$ such
that $x=\mathbf{1}\otimes b$.
\end{proof}

%%%%%%%%%%%%%%%%%%%%%%%%%%%%%%%%%%%%%%%%%%%%%%%%%%%%%%%%%%%%%%%%%%%%%%%%%%%%%%%%

\section{THE STRUCTURE OF CONSENSUS FUNCTION}
It is important to have  clear understanding of the structure of
consensus function in studying consensus problem.  Hence, in this
section,  we study the consensus function and present some
characterizations.

\subsection{Consensus Function is a Time-invariant Quantity}
We still consider   system (\ref{xt1}).  If it solves a consensus
problem for any initial state,  i.e.,  it satisfies the conditions
in Theorem \ref{dl21},  then for $\forall x \in \R^{mn}$, $\exists
b \in \R^m$ such that
$\mathop{\lim}\limits_{t\to\infty}e^{At}x=\mathbf{1}\otimes
b\triangleq x^*$,
  $x_i^*=[I_m, 0, \cdots, 0](\mathbf{1}\otimes b)=
[I_m, 0, \cdots, 0]\mathop{\lim}\limits_{t\to\infty}e^{At}x$,
where $I_m$ is identical matrix of order $m$.  Let $\chi(x)=[I_m$,
$0,  \cdots, 0] \mathop{\lim}\limits_{t\to\infty}e^{At}x$,  then
the system solves the $\chi$-consensus problem. It is easy to see
that the consensus function is determined by $A$.  Hence, if
system (\ref{xt1}) solves a consensus problem for any initial
state,  it must solve the $\chi$-consensus problem for some
consensus function $\chi$.

If $\mathop{\lim}\limits_{t\to\infty}e^{At}x=x^*$ is an
equilibrium for any $x \in \R^{mn}$,  we have $Ax^*=0$,  i.e.,
$A\mathop{\lim}\limits_{t\to\infty}e^{At}x=0$ for any $x \in
\R^{mn}$.  Thus $A\mathop{\lim}\limits_{t\to\infty}e^{At}=0$.
Since $e^{At}A=A e^{At}$, we have
$\mathop{\lim}\limits_{t\to\infty}e^{At}A=A\mathop{\lim}\limits_{t\to\infty}e^{At}=0$.
Hence $\frac{d\chi (x)}{dt}=[I,  0,  \cdots,
0]\mathop{\lim}\limits_{t\to \infty}e^{At}Ax= 0$.  So the
consensus function $\chi(x)$ is a time-invariant quantity. (Note
that $\mathop{\lim}\limits_{t\to \infty}e^{At}$ is a constant
matrix.)

\begin{rmrk}
If system (\ref{xt1}) solves the $\chi$-consensus problem for any
initial state, then the consensus problem can not be max- or min-
consensus. This is obvious by $\chi(x)=[I, 0, \cdots,
0]\mathop{\lim}\limits_{t\to \infty}e^{At}x$.
\end{rmrk}

\subsection{A Method to Obtain the Consensus Function}

For a given  system, the consensus function can be obtained  by
calculating $\mathop{\lim}\limits_{t\to \infty}e^{At}$.  When all
eigenvalues of $A$ have negative real parts, it is easy to obtain
that $\chi(x)\equiv 0$. However, if $0$ is an eigenvalue of $A$,
the calculation of $\mathop{\lim}\limits_{t\to \infty}e^{At}$
might be very complex. In what follows, we will illustrate that,
for some special cases, we can find a simple method to obtain the
consensus function.

Consider the following system:
\begin{equation}\label{xt2}
    \left[%
\begin{array}{c}
  \dot{x}_1 \\
  \dot{x}_2 \\
  \vdots \\
  \dot{x}_n \\
\end{array}%
\right]
=\left[%
\begin{array}{cccc}
  A_{11} & A_{12} & \cdots & A_{1n} \\
  A_{21} & A_{22} & \cdots & A_{2n} \\
  \vdots & \vdots & \ddots & \vdots \\
  A_{n1} & A_{n2} & \cdots & A_{nn} \\
\end{array}%
\right]
\left[%
\begin{array}{c}
  {x_1} \\
  {x_2} \\
  \vdots \\
  {x_n} \\
\end{array}
\right]
\end{equation}
denoted  by $\dot{x}=Ax$,  which satisfies the conditions in
Theorem \ref{dl21} and  rank$(A^2)=$rank$(A)=(n-1)m$.

By Theorem \ref{dl24}, it is easy to show that
\begin{equation}\label{a}
    \sum_{j=1}^n A_{ij}=0,  \forall i \in \n
\end{equation}

Let \[B=\mathop{\lim}\limits_{n\to \infty} e^{At}=\left[%
\begin{array}{cccc}
  B_{11} & B_{12} & \cdots & B_{1n} \\
  B_{21} & B_{22} & \cdots & B_{2n} \\
  \vdots & \vdots & \ddots & \vdots \\
  B_{n1} & B_{n2} & \cdots & B_{nn} \\
\end{array}%
\right]
=\left[%
\begin{array}{c}
  B_1 \\
  B_2 \\
  \vdots \\
  B_n \\
\end{array}%
\right],\] where $B_{ij}\in\R^{m\times m}, B_i\in\R^{m\times mn }$
for any $i, j \in \n$.

 Since system (\ref{xt2}) solves the $\chi$-consensus
problem for any $x(0) \in\R^{mn}$,  and let
\[x(0) =\underbrace{\left[%
\begin{array}{c}
  1 \\
  0 \\
  \vdots \\
  0 \\
\end{array}%
\right],  \left[\begin{array}{c}
  0 \\
  1 \\
  \vdots \\
  0 \\
\end{array}%
\right],  \cdots,  \left[\begin{array}{c}
  0 \\
  \vdots \\
  0 \\
  1 \\
\end{array}%
\right]}_{mn},\] respectively, we get $B_{11}=B_{21}= \cdots
=B_{n1}$, $B_{12}=B_{22}= \cdots =B_{n2}$,  $\cdots$,
$B_{1n}=B_{2n}= \cdots =B_{nn}$,  i.e.,  $B_{1}=B_{2}= \cdots
=B_{n}$.  We denote $B_i$ by $E=(E_1,  E_2,  \cdots,  E_n)$, where
$E_i\in\R^{m\times m}$ for any $i\in\n$, so
\[B=\left[%
\begin{array}{c}
  E\\
  E \\
  \vdots \\
  E \\
\end{array}%
\right]=\left[\begin{array}{cccc}
       E_1 & E_2 & \cdots & E_n \\
       E_1 & E_2 & \cdots & E_n \\
          \vdots & \vdots & \ddots & \vdots \\
       E_1 & E_2 & \cdots & E_n \\
        \end{array}\right].
\]

Since $\mathop{\lim}\limits_{t\to \infty}e^{At}A=A
\mathop{\lim}\limits_{t\to \infty}e^{At}=0_{(mn)\times (mn)}$,  we
get $BA=AB=0_{(mn)\times (mn)}$.  Therefore \[EA=0_{m\times (mn)}.
\]

Since $\chi(x)$ is an invariant quantity, we  have
\[\mathop{\lim}\limits_{t\to\infty}\chi(x(t))=\chi (\left[%
\begin{array}{c}
  Ex(0) \\
  Ex(0) \\
  \vdots \\
  Ex(0) \\
\end{array}%
\right] )=Ex(0).\] This implies
\[(E_1+E_2+\cdots+E_n)Ex(0)=Ex(0)
\] for all $x(0)\in \R^{mn}$.

By the theory of Jordan canonical form, we learn that rank$(B)=m$,
i.e., rank$(E)=m$,  so $\{ Ex(0)| \forall x(0) \in \R^{mn} \}=
\R^m$. Hence \begin{equation}\label{b} E_1+E_2+ \cdots +E_n=I_{m}.
\end{equation}

Because rank$(A)=(n-1)m$,  there exist $m$ linearly independent
vectors $\xi_1,  \xi_2,  \cdots,  \xi_m$ in $\R^{mn}$ such that
$\xi_i^T A=0$ for all $i\in\{1, 2, \cdots, m\}$. Let
\[Z=\left[%
\begin{array}{c}
  \xi_1^T \\
  \xi_2^T \\
  \vdots \\
  \xi_m^T \\
\end{array}%
\right]=[Z_1,  Z_2,  \cdots,  Z_n],\] where $Z_i \in \R^{m\times
m}$. Then there exists an invertible matrix $T\in \R^{m\times m}$
such that $E=TZ$.

By (\ref{b}), we have
 \[E\left[%
\begin{array}{c}
  I_m\\
  I_m \\
  \vdots \\
  I_m \\
\end{array}%
\right]=TZ\left[%
\begin{array}{c}
  I_m\\
  I_m \\
  \vdots \\
  I_m \\
\end{array}%
\right]=T(Z_1+Z_2+\cdots+Z_n)=I_m. \]
 So
$Z_1+Z_2+\cdots+Z_n$ is invertible and
$T=(Z_1+Z_2+\cdots+Z_n)^{-1}$.  Therefore
\[\chi(x)=Ex=(Z_1+Z_2+\cdots+Z_n)^{-1} Z x.\]

By the discussion above,  we get the following procedure to get
the consensus function:

\begin{ff}
\begin{enumerate}
    \item Choose arbitrarily $m$ linearly independent vectors $\xi_1,  \xi_2,  \cdots, \xi_m \in
    N(A^{T})$;
    \item Let $F=(\xi_1,  \xi_2,  \cdots, \xi_m)^{T}=(F_1,  F_2,  \cdots,
    F_m)$,  where $F_i\in \R^{m\times m}$, $i=1,  2,  \cdots,  m$;
    \item Let $T=F_1+F_2+\cdots+F_m$,  then $T$ is invertible;
    \item $\chi(x)=T^{-1}Fx$,  $\mathop{\lim}\limits_{n\to \infty}e^{At}=\left[%
\begin{array}{c}
  T^{-1}F \\
  T^{-1}F \\
  \vdots \\
  T^{-1}F \\
\end{array}%
\right].$
\end{enumerate}
\end{ff}

\subsection{An Example}

In what follows, we present an example  to show the effectiveness
of Method 1.

\begin{lz}\label{lz31}
\[A=\left[%
\begin{array}{cccc}
  -1 & 0 & 1 & 0 \\
  1 & 0 & -1 & 0 \\
  0 & 1 & 0 & -1 \\
  0 & 1 & 0 & -1\\
\end{array}%
\right].\]
\end{lz}

By Theorem \ref{dl21}, it is easy to verify that system
$\dot{x}=Ax$ solves a consensus problem for $m=2$.

Since $\left[%
\begin{array}{cccc}
  0 & 0 & -1 & 1 \\
  1 & 1 & 0 & 0 \\
\end{array}%
\right]A=0$,  let
\[T=\left[\begin{array}{cc}
  0 & 0  \\
  1 & 1  \\
\end{array}%
\right] + \left[\begin{array}{cc}
-1 & 1 \\
 0 & 0 \\
\end{array}%
\right]=\left[\begin{array}{cc}
 -1 & 1 \\
  1 & 1 \\
\end{array}%
\right],\]
 and thus \[T^{-1}=\left[\begin{array}{cc}
 -0. 5 & 0. 5 \\
  0. 5 & 0. 5 \\
\end{array}%
\right]. \]
Hence \[\chi(x)=T^{-1}\left[%
\begin{array}{cccc}
  0 & 0 & -1 & 1 \\
  1 & 1 & 0 & 0 \\
\end{array}%
\right]x\]
\[ =\left[%
\begin{array}{cccc}
  0. 5 & 0. 5 & 0. 5 & -0. 5 \\
  0. 5 & 0. 5 & -0. 5 & 0. 5 \\
\end{array}%
\right]x.\]

On the other hand, calculating $\mathop{\lim}\limits_{t\to\infty}
e^{At}$ directly,  we get
\[\mathop{\lim}\limits_{t\to\infty} e^{At}=\left[%
\begin{array}{cccc}
  0. 5 & 0. 5 & 0. 5 & -0. 5 \\
  0. 5 & 0. 5 & -0. 5 & 0. 5 \\
    0. 5 & 0. 5 & 0. 5 & -0. 5 \\
  0. 5 & 0. 5 & -0. 5 & 0. 5 \\
\end{array}%
\right].\] Hence the consensus function obtained by  Method 1 is
correct.

\subsection{Average Consensus Problem}

Average consensus problem has been discussed in \cite{c1} and
\cite{c2}, and the authors of them presented some necessary and
sufficient conditions. In this section, we also consider the
problem but from another viewpoint.

Based on the discussion in Subsection A,  we set average consensus
function
$\chi(x)=F_{m\times(mn)}x=\frac{1}{n}(x_1+x_2+\cdots+x_n)$  such
that $FA=0$, where $x\in\R^{mn}$, $x_i\in\R^m$ for any $i\in\n$.

Let
\[
x=\left[%
\begin{array}{c}
  1 \\
  0 \\
  \vdots \\
  0 \\
\end{array}%
\right],
x=\left[%
\begin{array}{c}
  0 \\
  1 \\
  \vdots \\
  0 \\
\end{array}%
\right], \cdots,
x=\left[%
\begin{array}{c}
  0 \\
  0 \\
  \vdots \\
  1 \\
\end{array}%
\right],
\]
respectively,  we get
\[
F=\left[%
\begin{array}{cccccccccccc}
  \frac{1}{n} &  &  &  & \frac{1}{n} &  &  &  &  &  &  &  \\
   & \frac{1}{n} &  &  &  & \frac{1}{n} &  &  & \cdot &\cdot  & \cdot &  \\
   &  & \ddots &  &  &  & \ddots &  &  &  &  &  \\
   &  &  & \frac{1}{n} &  &  &  & \frac{1}{n} &  &  &  &  \\
\end{array}%
\right. \]
\[ \left.
\begin{array}{cccccccccccc}
 \frac{1}{n} &  &  &  \\
 & \frac{1}{n} &  &  \\
 &  & \ddots &  \\
 &  &  & \frac{1}{n} \\
\end{array}%
\right].
\]
Therefore dim$(N(A^T))\geq m$.  By  Corollary \ref{tl22},  we have
the following lemma.

\begin{yl}\label{yl31}
If system (\ref{xt1}) solves the average consensus problem,  then
rank$(A)=(n-1)m$.
\end{yl}

\begin{dl}\label{dl34}
System (\ref{xt1}) solves the average consensus problem if and
only if rank$(A^2)$=rank$(A)=m(n-1)$,  $[I, I, \cdots, I]A=0$ and
$A[I, I, \cdots, I]^T=0$, where $I \in \R^{m\times m}$ is
identical matrix, and $0$ is an eigenvalue of $A$, and all the
other eigenvalues have negative real parts.
\end{dl}
\begin{proof}
The proof is obvious. We omit the details.
\end{proof}
\subsection{The Case of $m=1$}

In this subsection, we study the consensus problem in the case
$m=1$.

Consider the system
\begin{equation}\label{xt3}
\dot{x}=Ax,\end{equation}where  \[x=\left[%
\begin{array}{c}
  x_1 \\
  x_2 \\
  \vdots \\
  x_n \\
\end{array}%
\right]
, \quad A=\left[%
\begin{array}{cccc}
  a_{11}& a_{12}&\cdots&a_{1n} \\
  a_{11}& a_{12}&\cdots&a_{1n} \\
  \vdots &\vdots& \ddots&\vdots \\
  a_{11}& a_{12}&\cdots&a_{1n} \\
\end{array}%
\right],\]
 and $x_i,  a_{ij} \in \R,  \forall i,  j \in \n. $

\begin{dl}\label{dl35}
System (\ref{xt3}) solves the $\chi$-consensus problem for any
initial state if and only if\\
1) each eigenvalue of $A$ is $0$ or has negative real part;\\
2) if $0$ is an eigenvalue of $A$, then $A\mathbf{1}=0$, and {\rm
rank}$(A^2)$={\rm rank}$(A)=n-1$.
\end{dl}

So we get a general method to derive the consensus function.

\begin{ff}
Choose arbitrarily $y \in N(A^{T})$,  $y\neq 0$,  and let $y=[y_1,
y_2, \cdots, y_n]^T$,  then the consensus function is
\[\chi(x)=\frac{\sum_{i=1}^n y_i x_i}{\sum_{i=1}^n y_i},\]
where $x=[x_1,  x_2,  \cdots,  x_n]^T\in{R^n}$.
\end{ff}

We  have the following corollary.

\begin{tl}\label{tl36}
System (\ref{xt3}) solves the average consensus problem if and
only if all the nonzero eigenvalues of $A$ have negative real
parts, rank$(A^2)$=rank$(A)=n-1$, $\mathbf{1}^TA=0$, and
$A\mathbf{1}=0$. (\cite{c1},  Theorem 5)
\end{tl}

\begin{rmrk}
Naturally, if $A$ is a Laplacian matrix of some graph, then it
satisfies the conditions in Corollary \ref{tl36} except
$\mathbf{1}^TA=0$.
\end{rmrk}

\begin{rmrk}
we can view the consensus function
\[\chi(x)=\frac{\sum_{i=1}^n
y_i x_i}{\sum_{i=1}^n y_i}\]
 as a weighted average consensus
function.
\end{rmrk}

\section{CONSENSUS PROBLEM UNDER SIMILAR TRANSFORMATION}

Since consensus function is determined by $A$,  the study on the
structure of $A$ becomes an important issue.

Let
\[
J=\left[%
\begin{array}{ll}
0_{r\times r}&0_{r\times (mn-r)}\\
0_{(mn-r)\times r }&M_{(mn-r)\times (mn-r)}
\end{array}%
\right],
\]
where $M$ is a nonsingular real matrix,  the eigenvalues of $M$
have negative real parts and $r\leq m$.  Hence,  by Corollary
\ref{tl23}, if system (\ref{xt1}) solves the $\chi$-consensus
problem for any initial state, and dim$(N(A))=r$,  then there
exists a nonsingular real matrix $T$ such that $A=T^{-1}JT$.  We
will study the structure of $T$ in this section.

From the discussion in Section II, system (\ref{xt1}) solves the
$\chi$-consensus problem for any initial state if and only if
$\mathop{\lim}\limits_{t\to\infty}e^{At}$ exists and the
equilibriums have the form of $\mathbf{1}\otimes b$. In the
following, we will show that, for some systems, if
$\mathop{\lim}\limits_{n\to\infty}e^{At}$ exists,  we may find
similar transformation $T$ such that $T^{-1}AT$ solves the
$\chi$-consensus problem.

Consider the following system
\begin{equation}\label{xt41}
    \dot{x}=T^{-1}JTx.
\end{equation}

\begin{dl}\label{dl41}
System (\ref{xt41}) with $r=m$ solves the $\chi$-consensus problem
for any initial state if and only if
\[
T=\left[%
\begin{array}{cccc}
  T_{11} & T_{12} & \cdots & T_{1n} \\
  T_{21} & T_{22} & \cdots & T_{2n} \\
  \vdots & \vdots & \ddots & \vdots \\
  T_{n1} & T_{n2} & \cdots & T_{nn} \\
\end{array}%
\right],
\]
where $T_{ij}\in \R^{m\times m}$ for any $i, j \in \n$, satisfies
$T_U=\sum_{j=1}^nT_{1j}$ is invertible, $\sum_{j=1}^nT_{ij}=0$ for
$i=2,  \cdots,  n$,  and
\[
T_D=\left[%
\begin{array}{ccc}
 T_{22} & \cdots & T_{2n} \\
 \vdots & \ddots & \vdots \\
T_{n2} & \cdots & T_{nn} \\
\end{array}%
\right]
\]
is invertible.
\end{dl}

\begin{proof}
{\it Sufficiency}.  It suffices to prove that all equilibriums of
$T^{-1}JT$ have the form: $\mathbf{1}\otimes b$.

Let $b_1,  b_2,  \cdots,  b_m \in \R^m$ be linearly independent
vectors, and let $\zeta_1=\mathbf{1}\otimes b_1,
\zeta_2=\mathbf{1}\otimes b_2, \cdots,  \zeta_m=\mathbf{1}\otimes
b_m$, which are also linearly independent,  then
$T^{-1}JT\zeta_i=0$, $i=1, 2, \cdots, m$.  The sufficiency is
proved.

{\it Necessity}.  The equation $T^{-1}JTx=0$ must have $m$
linearly independent solutions:
\[
\zeta_1=\mathbf{1}\otimes b_1,  \zeta_2=\mathbf{1}\otimes b_2,
\cdots, \zeta_m=\mathbf{1}\otimes b_m,
\] where $b_i\in{\R^m}$, $i=1, \cdots, m.$

Notice that
\[T^{-1}JTx=0\Leftrightarrow JTx=0 \Leftrightarrow Tx=\left[%
\begin{array}{c}
  c \\
  0 \\
  \vdots \\
  0\\
\end{array}%
\right]\]
 for some $c \in \R^{m}$.  So there exist $c_1,  c_2,  \cdots,  c_m \in
 \R^m$ such that
\[
T\zeta_1=\left[%
\begin{array}{c}
  c_1 \\
  0 \\
  \vdots \\
  0\\
\end{array}%
\right],
T\zeta_2=\left[%
\begin{array}{c}
  c_2 \\
  0 \\
  \vdots \\
  0\\
\end{array}%
\right], \cdots,
T\zeta_m=\left[%
\begin{array}{c}
  c_m \\
  0 \\
  \vdots \\
  0\\
\end{array}%
\right],
\]
and
\[
\left\{%
\begin{array}{ll}
  \big(\sum\limits_{j=1}^n T_{1j}\big) b_k=c_k, &  \\
  \big(\sum\limits_{j=1}^n T_{ij}\big) b_k=0, & i=2, 3, \cdots, n, \\
\end{array}%
\right.
\]
for $k=1, 2, \cdots, m$. Since $b_1, b_2, \cdots, b_m$ are
linearly independent,  we have
\[
\sum\limits_{j=1}^n T_{ij}=0,  i=2, 3, \cdots, n.
\]
Since $T$ is invertible,  we obtain that $T_U$ and $T_D$ are
invertible.
\end{proof}

For the average consensus problem,  we have the following theorem.
\begin{dl}\label{dl42}
System (\ref{xt41}) solves the average consensus problem for any
initial state if and only if $r=m$ and
\[
T=\left[%
\begin{array}{cccc}
  T_{11} & T_{11} & \cdots & T_{11} \\
  T_{21} & T_{22} & \cdots & T_{2n} \\
  \vdots & \vdots & \ddots & \vdots \\
  T_{n1} & T_{n2} & \cdots & T_{nn} \\
\end{array}%
\right],
\]
where $\sum_{j=1}^{n}T_{ij}=0$,  $i=2, 3, \cdots, n$,  and
$T_{11}$ and $T_D$ all are invertible.
\end{dl}
\begin{proof}
{\it Sufficiency}.  Let
\[T^{-1}=S=\left[%
\begin{array}{cccc}
  S_{11} & S_{12} & \cdots & S_{1n} \\
  S_{21} & S_{22} & \cdots & S_{2n} \\
  \vdots & \vdots & \ddots & \vdots \\
  S_{n1} & S_{n2} & \cdots & S_{nn} \\
\end{array}%
\right],\] where $S_{ij}\in\R^{m\times m}$ for any $i, j\in \n$.
We have
\[S_{11}=S_{21}=\cdots S_{n1}=\frac{1}{n}T_{11}^{-1},\]
\[\sum_{i=1}^nS_{ij}=0,  j=2, 3,  \cdots,  n,\]
and
\[
S_D=\left[%
\begin{array}{ccc}

 S_{22} & \cdots & S_{2n} \\
 \vdots & \ddots & \vdots \\
 S_{n2} & \cdots & S_{nn} \\
\end{array}%
\right]
\]
is invertible.

It is obvious that \[T^{-1}JT\left[%
\begin{array}{c}
  I \\
  I \\
  \vdots \\
  I \\
\end{array}%
\right]=0\]
 and
 \[
 [I,  I,  \cdots,  I]T^{-1}JT=0.
 \]
By Theorem \ref{dl34},  $T^{-1}JT$ solves the average consensus
problem.

{\it Necessity}.  If system (\ref{xt41}) solves the average
consensus problem,  we have
\[
T^{-1}\lim_{t\to\infty}e^{Jt}T=\frac{1}{n}\left[%
\begin{array}{cccc}
  I & I & \cdots & I \\
  I & I & \cdots & I \\
  \vdots & \vdots & \ddots & \vdots \\
  I & I & \cdots & I \\
\end{array}%
\right],
\]
which implies \[ S\left[%
\begin{array}{cccc}
  I & 0 & \cdots & 0 \\
  0 & 0 & \cdots & 0 \\
  \vdots & \vdots & \ddots & \vdots \\
  0 & 0 & \cdots & 0 \\
\end{array}%
\right] T\]
\[=\left[%
\begin{array}{cccc}
  S_{11}T_{11} & S_{11}T_{12} & \cdots & S_{11}T_{1n} \\
  S_{21}T_{11} & S_{21}T_{12} & \cdots & S_{21}T_{1n} \\
  \vdots & \vdots & \ddots & \vdots \\
 S_{n1}T_{11} &S_{n1}T_{12} & \cdots &S_{n1}T_{1n} \\
\end{array}%
\right]
\]
\[=\frac{1}{n}\left[%
\begin{array}{cccc}
  I & I & \cdots & I \\
  I & I & \cdots & I \\
  \vdots & \vdots & \ddots & \vdots \\
  I & I & \cdots & I \\
\end{array}%
\right].\] We obtain that $T_{11}=T_{12}=\cdots=T_{1n}$ are
invertible.

By Lemma \ref{yl31},  we get $r=m$, and by Theorem \ref{dl41},  we
have $\sum_{j=1}^nT_{ij}=0$,  $i=2, 3, \cdots, n$ and $T_D$ is
invertible.
\end{proof}

For the case of $r<m$,  $T$ is very complex,  but we still have
the following theorem.

\begin{dl}
If
\[
T=\left[%
\begin{array}{cccc}
  T_{11} & T_{12} & \cdots & T_{1n} \\
  T_{21} & T_{22} & \cdots & T_{2n} \\
  \vdots & \vdots & \ddots & \vdots \\
  T_{n1} & T_{n2} & \cdots & T_{nn} \\
\end{array}%
\right],
\]
where $T_U$ and $T_D$ defined as in Theorem 4.1 are  invertible,
and $\sum_{j=1}^nT_{ij}=0$, $i=2, 3, \cdots, n$,  then  system
(\ref{xt41}) solves the $\chi$-consensus problem for any $r\leq
m$.
\end{dl}
\begin{proof}
For arbitrary $r$ linearly independent vectors $b_1,  b_2, \cdots,
b_r\in \R^r$, we define
\[
\zeta_i=\mathbf{1}\otimes\left(T_U^{-1}\left[%
\begin{array}{c}
 b_i \\
  0 \\
\end{array}%
\right]_{m\times 1}\right)\] for $i=1, 2, \cdots, r.$ Then
$\zeta_1, \zeta_2,  \cdots,  \zeta_r$ are linearly independent and
$T^{-1}JT\zeta_i=0$, $i=1, 2, \cdots, r$.

Therefore, system (\ref{xt41}) solves the $\chi$-consensus
problem.
\end{proof}

\begin{rmrk}
If $\lim_{t\to\infty}e^{At}$ exists and dim$(N(A))=r\leq m$,  then
we can find an invertible matrix $T$  such that system
$\dot{x}=T^{-1}ATx$ solves the $\chi$-consensus problem.

Moreover, if $r=m$,  we can  find $T$ such that system
$\dot{x}=T^{-1}ATx$ solves the average consensus problem.
\end{rmrk}

\begin{lz}
For any initial state, the system in Example \ref{lz31} solves a
consensus problem, but not the average consensus problem. We will
provide the procedure to find invertible $T$ such that system
$\dot{x}=T^{-1}ATx$ solves the average consensus problem.

First, we can choose an invertible matrix $T_1$,
\[T_1= \left [\begin {array}{cccc}
2&0&1&-1\\
6&4&-1&1\\
2&0&-1&-1\\
6&4&-1&-1
\end{array}\right ],\]
 such that
 \[
 B=T_1^{-1}AT_1=\left [\begin {array}{cccc} 0&0&0&0\\
 0&0&0&0\\
 0&0&-1&-1\\
0&0&1&-1\end {array}\right],
 \]
 which has the form of $J$.
\end{lz}

By Theorem \ref{dl42}, we let
 \[ T_2=\left[\begin {array}{cccc} 1& 0& 1& 0\\
0&1&0 &1\\
-1&0&1&0\\
0&-1&0&1\end {array}\right].
 \]
Then
\[
C=T_2^{-1}BT_2=\left[\begin {array}{cccc}
-0.5&-0.5&0.5&0.5\\
0.5&-0.5&-0.5&0.5\\
0.5&0.5&-0.5&-0.5\\
-0.5&0.5&0.5&-0.5
\end{array}\right].
\]

By Theorem \ref{dl34}, $\dot{x}=Cx$ solves the average consensus
problem for any initial state. Let \[T=T_1T_2=\left [\begin
{array}{cccc}
1&1&3&-1\\
7&3&5&5\\
3&1&1&-1\\
7&5&5&3\end {array}\right ],\] then $\dot{x}=T^{-1}ATx$ solves the
average consensus problem for any initial state.

\section{CONSENSUS PROBLEM UNDER ARBITRARY SWITCHING}

In this section, we investigate the consensus problem of system
\begin{equation}\label{xt51}
\dot{x}=A(t)x.
\end{equation}
The study on the consensus problem of (\ref{xt51}) is difficult.
Here, we only consider some special cases. We view (\ref{xt51}) as
a switched system and $A(t)$ is a constant matrix in each
switching interval.

We consider the following system, each subsystem of which is the
same as (\ref{xt1}),
\begin{equation}\label{xt4}
    \dot{x}(t)=A_{s(t)}x(t),
\end{equation}
where $s(t): \R^+ \to \n$ is the switching signal.

Generally speaking,  not all switched systems    solve a consensus
problem for any initial state. But  some special switched systems
can solve a consensus problem.

We assume that, for $\forall s \in \n$,

1) $\dot{x}=A_sx$ solves a consensus problem;

2) $A_s=\left[%
\begin{array}{cccc}
  A_{s11} & A_{s12} & \cdots & A_{s1n} \\
  A_{s21} & A_{s22} & \cdots & A_{s2n} \\
  \vdots & \vdots & \ddots & \vdots \\
  A_{sn1} & A_{sn2} & \cdots & A_{snn} \\
\end{array}%
\right],$\\
 where $A_{sij}$ is a symmetric and positive definite matrix
(denoted by $A_{sij}>0$) for all $i,  j \in \n,  i \neq  j$;

3) $A_{sii}=-\sum\limits_{j=1, j\neq i}^{n}A_{sij} $;

4) every subsystem  has the same consensus function
$\chi(x)=F_{m\times mn}x=[F_1,  F_2,  \cdots,  F_n]x$,  where
$F_i\in \R^{m\times m}$, $F_i>0$ for all $i \in \n$, and $FA_s=0$;

5) $F_iA_{sij}=A_{sij}F_i$ for all $i, j\in\n$, which implies
$F_iA_{sij}>0$.

Before presenting Theorem \ref{dl52}, we first prove the following
lemma.

\begin{yl}
Let
\begin{equation}\label{eq51}
    L=\left[%
\begin{array}{cccc}
  L_{11} & -L_{12} & \cdots & -L_{1n} \\
  -L_{21} & L_{22} & \cdots & -L_{2n} \\
  \vdots & \vdots & \ddots & \vdots \\
  -L_{n1} & -L_{n2} & \cdots & L_{nn} \\
\end{array}%
\right]
\end{equation}
be a symmetric matrix,  where $L_{ij}\in \mathbf{R}^{m\times m}$,
$L_{ij}>0$ and $\sum_{j=1, j\neq i}^n L_{ij}=L_{ii}$ for all $i,
j\in \n$.

If the eigenvalues of $L$ are arranged in an increasing order
$\lambda_1 \leq \lambda_2 \leq \cdots \leq \lambda_{mn}$,  then we
have
\[\lambda_1 = \lambda_2=\cdots=\lambda_{m}=0,  \lambda_{m+1}>0.\]
Then we call $L$  a {\it block laplacian matrix}.
\end{yl}

\begin{proof}
For any $x\in \R^{mn}$, we have
\[
x^TLx=\sum_{i=1}^n x_i^TL_{ii}x_i-\sum_{i\neq j}x_i^TL_{ij}x_j\]
\[
=\sum_{i\neq j}(x_i^TL_{ij}x_i)-\sum_{i\neq j}x_i^TL_{ij}x_j =
\sum_{i\neq j}(x_i^TL_{ij}x_i-x_i^TL_{ij}x_j)
\]
\[
=\frac{1}{2}\sum_{i\neq
j}(x_i^TL_{ij}x_i+x_j^TL_{ij}x_j-x_i^TL_{ij}x_j-x_j^TL_{ij}x_i)
\]
\[=\frac{1}{2}\sum_{i\neq j}(x_i-x_j)^TL_{ij}(x_i-x_j)\]
\[=\frac{1}{2}\sum_{i,  j}(x_i-x_j)^TL_{ij}(x_i-x_j).\]

Therefore $\lambda_1 = \lambda_2=\cdots=\lambda_{m}=0,
\lambda_{m+1}>0$.
\end{proof}

\begin{dl}\label{dl52}
If  switched system  (\ref{xt4}) satisfies the conditions
(1)--(5), then, under arbitrary switching, the solution of the
system globally asymptotically converges to
$\mathbf{1}\otimes\chi(x(0))$,  i.e., the switched system solves
the $\chi$-consensus problem.
\end{dl}
\begin{proof} Since $\chi(x)$ is an invariant
quantity for every subsystem,  $\chi(x)$ is also an invariant
quantity under switching.  For any solution $x(t)$,  let
$x(t)=\mathbf{1}\otimes \chi(x)+ \delta(t)$.  We refer to $\delta$
as the {\it  (group) disagreement vector}. Then
\[Fx=F(\mathbf{1}\otimes \chi(x))+F\delta = Fx+ F\delta.\] Hence,
we have
\begin{equation}\label{del}
    F\delta=0.
\end{equation}

For a given $s\in \n$, since $\dot{x}=A_sx$,  we have
\[\mathbf{1}\otimes\frac{d\chi(x)}{dt}+\dot{\delta} =
A_s(\mathbf{1}\otimes\chi(x))+A_s{\delta},\]
which implies
\[\dot{\delta}=A_s\delta.\]

Let \[\Theta=\left[%
\begin{array}{cccc}
 F_1 &  &  &  \\
   & F_2 &  &  \\
   &  & \ddots &  \\
   &  &  & F_n \\
\end{array}%
\right],\] then $\Theta>0$.

 Let $V(\delta)=\delta^{T}\Theta \delta$,  then
  \[\frac{dV}{dt}=2\delta^{T} \Theta \dot{\delta}= \delta^{T}( \Theta
  A_s+A_s^{T}\Theta^{T})\delta.\]

  Let $-L=\Theta  A_s+A_s^{T}\Theta^{T}$,  then $L$ is a block laplacian
  matrix
  by assumption.  We can easily get that
  $N(L)=R(\mathbf{1}\otimes I)$ and $N(L) \cap N(F)=R(\mathbf{1}\otimes I) \cap
  N(F)=\{0\}$.  We  divide the linear space $\R^{mn}$ into the direct
  sum of
  $R(\mathbf{1}\otimes I)$ and its orthogonal complement space
  $R(\mathbf{1}\otimes I)^{\bot}$,  then we have
  \[ \R^{mn}=R(\mathbf{1}\otimes I)\oplus R(\mathbf{1}\otimes I)^{\bot}.\]
  Correspondingly,  $\delta=\delta_1 +\delta_2$,  $\delta_1 \in R(\mathbf{1}\otimes I),  \delta_2 \in
  R(\mathbf{1}\otimes I)^{\bot}$.  Let $P$ be the orthogonal projector from
  $\R^{mn}$ onto $R(\mathbf{1}\otimes I)^{\bot}$ such that $\delta_2=P\delta$.
  Since $F\delta=0$,  we have  $P\delta=\delta_2\neq 0$  if $\delta\neq 0$.

  Hence $\frac{dV}{dt}= -\delta^T L \delta=
   -(\delta_1+\delta_2)^T L(\delta_1+\delta_2)$
    $=-\delta_2^T P^TLP \delta_2$ $\leq -\lambda_{m+1} \delta_2^T P^T
  P \delta_2$
  $<0$,  where $\lambda_{m+1}>0$ is the $(m+1)$th smallest eigenvalue
of
  $L$.  This shows that $V(\delta(t))$ is a valid common Lyapunov
  function for the group-disagreement,  i.e. ,  under arbitrary switching, the
switched system solves the $\chi$-consensus problem.
\end{proof}

\begin{rmrk}
The assumptions (1-5) seem rather strict,  but this kind of system
really exists extensively.  For example,  it is easy to show that
the system
\begin{equation}\label{rm53}
    \dot{x}=-Lx,
\end{equation}
where $L$ is block laplacian matrix,  solves the average consensus
problem and satisfies the assumptions (1--5).
\end{rmrk}

\section{CONCLUSIONS}

For linear dynamic systems,  consensus problem has been discussed
from a new viewpoint.  The structure of the consensus functions
has been characterized.  An example has been presented to
illustrate the effectiveness of our results.  Some necessary and
sufficient conditions for consensus problem under similar
transformation have also been obtained.  Finally,  we characterize
a class of dynamic switched systems that  solve  a consensus
problem under arbitrary switching.

\section{APPENDIX: GRAPH THEORY PRELIMINARIES}

In this section,  we briefly summarize some basic concepts and
results in graph theory that are useful in dealing with the
consensus problem.  More comprehensive discussions can be found in
\cite{C.  Godsil and G.  Royle}.

A {\it undirected graph} $\cal{G}$ consists of a {\it vertex set}
${\cal{V}}=\{n_1,  n_2,  \cdots,  n_m\}$ and an {\it edge set}
${\cal{E}}=\{(n_i,  n_j): n_i,  n_j \in {\cal{V}}\}$,  where an
{\it edge} is an unordered pair of distinct vertices of $\cal{V}$.
If $n_i,  n_j\in{\cal{V}}$,  and $(n_i,  n_j)\in{\cal{E}}$,  then
we say that $n_i$ and $n_j$ are {\it adjacent} or {\it neighbors}.
An {\it oriented graph} is a graph together with a particular
orientation,  where the {\it orientation} of a graph ${\cal {G}}$
is the assignment of a direction to each edge,  so edge $(n_i,
n_j)$ is an directed edge (arc) from $n_i$ to $n_j$.  The {\it
incidence matrix} $B$ of an oriented graph $\cal{G}$ is the $\{0,
\pm 1\}$-matrix with rows and columns indexed by the vertices and
edges of $\cal{G}$,  respectively,  such that the $ij$-entry is
equal to 1 if edge $j$ is ending on vertex $n_i$,  -1 if edge $j$
is beginning with vertex $n_i$,  and 0 otherwise.  Define the {\it
Laplacian matrix} of $\cal{G}$ as $L{(\cal{G})}=BB^T. $
$L{(\cal{G})}$ is always positive semi-definite.  Moreover,  for a
connected graph,  $L{(\cal{G})}$ has a single zero eigenvalue, and
the associated right eigenvector is ${\bf 1}_m$.

\end{document}